\documentclass{amsart}
\usepackage{amscd,amsmath,amssymb,amsthm,latexsym,pinlabel}

\def\co{\colon\thinspace}
\newcommand{\C}{{\mathbb C}}

\newcommand{\N }{{\mathbb N}}
\newcommand{\R}{{\mathbb R}}
\newcommand{\Z}{{\mathbb Z}}
\newcommand{\Hy}{{\mathbb H}}

\newcommand{\bw}{\mathbf{w}}

\newcommand{\orb}{\mbox{\scriptsize\rm orb}}
\newcommand{\ou}{\overline{u}}
\newcommand{\ov}{\overline{v}}
\newcommand{\oq}{\overline{q}}
\newcommand{\orho}{\overline{\rho}}
\newcommand{\otheta}{\overline{\vartheta}}
\newcommand{\PSL}{\mbox{\rm PSL}_2\R}
\newcommand{\Hom}{\mbox{\rm Hom}}
\newcommand{\SL}{\widetilde{\mbox{\rm SL}}_2}
\newcommand{\calR}{{\mathcal R}}
\newcommand{\orh}{\overline{\rho}}

\newcommand{\tf}{\tilde{f}}
\newcommand{\tg}{\tilde{g}}
\newcommand{\tu}{\tilde{u}}
\newcommand{\tv}{\tilde{v}}
\newcommand{\tq}{\tilde{q}}
\newcommand{\tih}{\tilde{h}}
\newcommand{\tpi}{\widetilde{\pi}}
\newcommand{\trho}{\widetilde{\rho}}
\newcommand{\ttheta}{\widetilde{\vartheta}}
\newcommand{\tSigma}{\widetilde{\Sigma}}
\newcommand{\jth}{j^{\mbox{\scriptsize\rm th}}}
\newcommand{\rth}{r^{\mbox{\scriptsize\rm th}}}
\newcommand{\Out}{\mbox{\rm Out}}
\newcommand{\Aut}{\mbox{\rm Aut}}
\newcommand{\Inn}{\mbox{\rm Inn}}
\newcommand{\T}{{\mathcal T}}
\newcommand{\M}{{\mathcal M}}

\newtheorem{thm}{Theorem}
\newtheorem{lem}[thm]{Lemma}
\newtheorem{prop}[thm]{Proposition}

\theoremstyle{definition}
\newtheorem*{defn}{Definition}
\newtheorem*{rem}{Remark}
\newtheorem*{rems}{Remarks}

\newtheorem*{claim}{Claim}

\begin{document}

\title[Generalised spin structures on orbifolds]{Generalised spin
structures on\\2-dimensional orbifolds}

\author{Hansj\"org Geiges}
\address{Mathematisches Institut, Universit\"at zu K\"oln,
Weyertal 86--90, 50931 K\"oln, Germany}
\email{geiges@math.uni-koeln.de}

\author[Jes\'us Gonzalo]{Jes\'us Gonzalo P\'erez}
\address{Departamento de Matem\'aticas,
Universidad Aut\'onoma de Ma\-drid, 28049 Madrid, Spain.}
\email{jesus.gonzalo@uam.es}
\thanks{J.~G. is partially supported by
         grants MTM2007-61982 from MEC Spain and MTM2008-02686 from MICINN
Spain.}

\date{}

\begin{abstract}
Generalised spin structures, or $r$-spin structures,
on a $2$-dimen\-sio\-nal orbifold $\Sigma$ are $r$-fold fibrewise connected
coverings (also called $\rth$ roots)
of its unit tangent bundle $ST\Sigma$. We investigate such
structures on hyperbolic orbifolds. The conditions on $r$ for such
structures to exist are given. The action of the diffeomorphism
group of $\Sigma$ on the set of $r$-spin structures is described,
and we determine the number of orbits under this action and their size.
These results are then applied to describe the moduli space of taut
contact circles on left-quotients of the $3$-dimensional
geometry $\SL$.
\end{abstract}



\maketitle

\section{Introduction}
Spin structures on manifolds have been studied extensively, not least
because of their relevance to physics. A spin structure on
a Riemann surface $\Sigma$ may be thought of as a square root of the tangent
bundle $T\Sigma$, that is, a holomorphic line bundle $\mathcal{L}$
with $\mathcal{L}\otimes\mathcal{L}=T\Sigma$. On the level of
the unit tangent bundle $ST\Sigma$, a spin structure can be
interpreted as a fibrewise connected double covering
$M\rightarrow ST\Sigma$ by another $S^1$-bundle $M$ over~$\Sigma$.

It is this last definition which most easily generalises to
2-dimensional orbifolds and coverings of higher order. This
is not just generalisation for generalisation's sake.
For instance, such objects appear in the work of Witten~\cite{witt93}
on matrix models of $2$-dimensional quantum gravity, see
also~\cite{poli04}. Here the viewpoint is that of Algebraic
Geometry, where an $\rth$ root of the tangent bundle
of a Riemann surface $\Sigma$ is considered to be a holomorphic
line bundle whose $\rth$ tensor power equals $T\Sigma$.
In that framework, questions of moduli have been
studied by Jarvis~\cite{jarv00} and others.

Our personal motivation for investigating such $\rth$ roots
comes from the moduli problem for taut contact circles
on $3$-manifolds. These structures were introduced in~\cite{gego95a},
where we also classified the $3$-manifolds which admit such structures.
The moduli question was largely settled in~\cite{gego02},
but certain details as to the precise geometry of the moduli spaces
had been left open. These details hinge on the classification of
$\rth$ roots of the unit tangent bundle of $2$-dimensional
hyperbolic orbifolds.

Here is an outline of the paper. In Section~\ref{section:orbifold}
we present the basics of $2$-dimensional hyperbolic orbifolds,
mostly to set up notation. In Section~\ref{section:tangent}
we recall the definition of the unit tangent bundle of an
orbifold. Roots of such unit tangent bundles are defined
in Section~\ref{section:roots}, where we determine the conditions
on $r$ (in terms of the genus and multiplicities
of the cone points of the orbifold~$\Sigma$) for $\rth$ roots to exist.
We also set up a one-to-one correspondence between $\rth$ roots
and certain homomorphisms on the fundamental group of
the unit tangent bundle of~$\Sigma$
(Theorem~\ref{thm:existence}). In Section~\ref{section:diffeo}
this is used to investigate the action by the diffeomorphism group
of $\Sigma$ on the set of $\rth$ roots. The number of orbits under
this action is determined (Proposition~\ref{prop:standard}),
as well as the length of the orbits (Proposition~\ref{prop:number}).
In Section~\ref{section:algebraic} we reformulate this
action by the diffeomorphism group in algebraic terms
as an action by the outer automorphism group of the orbifold
fundamental group. Finally, in Section~\ref{section:moduli}
we use this algebraic reformulation and the results of the previous
sections to describe the Teichm\"uller space (Theorem~\ref{thm:teich})
and moduli space (Theorem~\ref{thm:moduli}) of taut contact circles
on left-quotients of the $3$-dimensional Thurston geometry~$\SL$.
In particular, we are interested in the enumeration of the
connected components of the moduli space; this gives the number
of distinct taut contact circles up to diffeomorphism and
deformation.

Sections \ref{section:orbifold} to \ref{section:diffeo}
are completely self-contained. The final two sections depend
to some degree on our earlier work~\cite{gego02}, but except for
the algebraic reformulation of the moduli problem we do not
need to quote any details from that paper.
\section{Hyperbolic orbifolds}
\label{section:orbifold}
Throughout this paper, let $\Sigma$ be a fixed (closed, orientable,
$2$-dimensional) orbifold
of genus $g$ and with $n$ cone points of multiplicity $\alpha_1,\ldots ,
\alpha_n$. Moreover, it is assumed that $\Sigma$ is of
hyperbolic type, i.e.\ its orbifold Euler characteristic, defined as
\[ \chi^{\orb}(\Sigma )=2-2g-n+\sum_{j=1}^n\frac{1}{\alpha_j},\]
is assumed to be negative. This condition on the orbifold
Euler characteristic determines those orbifolds which admit
a hyperbolic metric; however, as yet we do not fix such
a hyperbolic structure.

\begin{figure}[h]
\labellist
\small\hair 2pt
\pinlabel $\oq_n$ [br] at 470 134
\pinlabel $\oq_1$ [bl] at 141 134
\pinlabel $\widetilde{x}_0$ [l] at 600 187
\pinlabel $\ou_1$ [tr] at 512 296
\pinlabel $\ov_g$ [tl] at 102 296
\pinlabel $\ov_1$ [tr] at 449 347
\pinlabel $\ou_g$ [tl] at 165 347
\endlabellist
\centering
\includegraphics[scale=0.4]{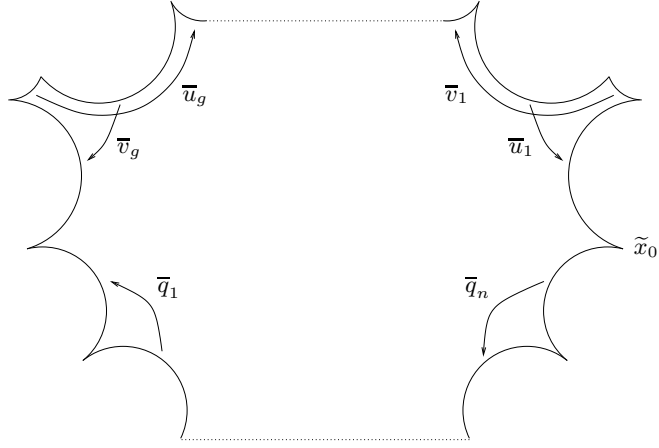}
  \caption{A fundamental domain for $\Sigma$.}
  \label{figure:domain1}
\end{figure}

The orbifold fundamental group $\pi^{\orb}$ of $\Sigma$ is
defined as the deck transformation group of the universal
covering $\tSigma\rightarrow\Sigma$. We briefly
recall the geometric realisation of this group and its
standard presentation. To that end,
choose a base point $x_0\in\Sigma$ distinct from all the cone points,
and a lift $\widetilde{x}_0\in\tSigma$ of $x_0$
in the universal covering space. Choose a system
of $2g$ loops on~$\Sigma$, based at~$x_0$, and a curve from $x_0$ to
each of the cone points, such that $\Sigma$ looks as in
Figure~\ref{figure:domain1}
when cut open along these $2g+n$ curves. We may interpret that figure
as a fundamental region in~$\tSigma$; it is determined
(amongst all possible fundamental regions whose boundary polygon
maps to the chosen system of curves) by the indicated placement of
$\widetilde{x}_0$ on its boundary. Notice that the sides
of this polygon identified by the deck transformation $\oq_j$ meet
at a vertex mapping to the $\jth$ cone point in~$\Sigma$; all other
vertices are lifts of~$x_0$.

Let $\ou_1,\ov_1,\ldots ,\ou_g,\ov_g,\oq_1,\ldots ,\oq_n$ be the
deck transformations of $\tSigma$ which effect the gluing
maps of the sides of the chosen fundamental polygon as indicated
in Figure~\ref{figure:domain1}. From the figure we see that
the deck transformation
$\prod_i[\ou_i,\ov_i]\prod_j \oq_j$
(read from the right as a composition of maps) fixes the
point~$\widetilde{x}_0$, which is not the lift of a cone point,
so we conclude
\[ \prod_i[\ou_i,\ov_i]\prod_j \oq_j=1.\]
Similarly, we have
\[ \oq_j^{\alpha_j}=1,\;\; j=1,\ldots ,n.\]
These relations give the standard presentation of $\pi^{\orb}$ as
\[ \pi^{\orb}  =  \Bigl\{ \ou_1,\ov_1,\ldots ,\ou_g,\ov_g,\oq_1,\ldots ,
\oq_n\co \prod_i[\ou_i,\ov_i]\prod_j\oq_j=1,\,\oq_j^{\alpha_j}=1\Bigr\} .\]

Once $\Sigma$ has been
equipped with a hyperbolic structure and an orientation,
then $\tSigma ={\mathbb H}^2$ and the $\ou_i,\ov_i,\oq_j$
are orientation preserving isometries of~${\mathbb H}^2$, i.e.\ elements
of $\PSL$. The identification of $\tSigma$ with
$\Hy^2$ is uniquely determined if we specify, for instance,
the lift $\widetilde{x}_0\in\Hy^2$, the initial direction
of one of the edges of the fundamental polygon
emanating from that point, and require that the orientation
lifted from $\Sigma$ coincide with a chosen orientation
of~$\Hy^2$. In this way an oriented hyperbolic structure
on $\Sigma$ defines an element of the Weil space
$\calR (\pi^{\orb},\PSL )$ of faithful representations of $\pi^{\orb}$
in $\PSL$ with discrete and cocompact image. Conversely,
any representation $\orh\in\calR (\pi^{\orb},\PSL )$ determines
a diffeomorphic copy $\orh (\pi^{\orb})\backslash\Hy^2$
of $\Sigma$ with a hyperbolic structure and an orientation.

It is possible to designate one of the orientations on any
given $\Sigma$ as {\em positive\/} and the other as
{\em negative\/} in the following way. If there are cone points, it suffices
to observe that $\orho (\oq_j)$ is a rotation by
$\pm 2\pi/\alpha_j$, with the same sign for each $j=1,\ldots ,n$.
(The sign is well defined even for $\alpha_j=2$ when we regard the rotation
as being through the interior of the fundamental domain.)
Observe in Figure~\ref{figure:domain1} how the direction of
rotation around the cone points relates to the orientation
given by the pairs of arrows indicating the action of $\ou_i$ and $\ov_i$;
thus, any such pair of arrows allows us to determine the orientation of
$\Sigma$, also when no cone points are present.
We write $\calR^{\pm}(\pi^{\orb},\PSL )$ for the
corresponding components of $\calR (\pi^{\orb},\PSL )$.
Any two representations $\orho_1,\orho_2\in \calR (\pi^{\orb},\PSL )$
are related via conjugation with a diffeomorphism of~$\Hy^2$.
This diffeomorphism will be orientation preserving or reversing,
depending on whether $\orho_1,\orho_2$ lie in the same component
or not, see also~\cite[pp.~59/60]{gego02}. This orientation
issue will only become relevant in Section~\ref{section:moduli}
of the present paper.
\section{The unit tangent bundle of an orbifold}
\label{section:tangent}
The unit tangent bundle of an oriented hyperbolic orbifold $\Sigma$ is
defined as follows, see~\cite[p.~466]{scot83}. Write
$\SL$ for the universal cover of $\PSL$. There is
a short exact sequence
\[ 0\longrightarrow\Z\longrightarrow\SL
\stackrel{p}{\longrightarrow}\PSL\longrightarrow 1.\]
Realise the given hyperbolic structure and orientation on $\Sigma$ by
a choice of representation $\orh\in\calR (\pi^{\orb},\PSL )$.
Then set
\[ ST\Sigma =p^{-1}(\orh (\pi^{\orb}))\backslash\SL ;\]
this is the unit tangent bundle of~$\Sigma$. It is in
a natural way the total space of a Seifert bundle over~$\Sigma$
with normalised Seifert invariants
\[ \{ g;b=2g-2; (\alpha_1,\alpha_1-1),\ldots ,(\alpha_n,\alpha_n-1)\} .\]

\begin{rem}
There is a tricky orientation issue here. The group $\PSL$ of
orientation preserving isometries of $\Hy^2$ acts,
via the differential,
transitively and with trivial point stabilisers on the unit tangent bundle
$ST\Hy^2$ of $\Hy^2$ (see Scott's survey~\cite{scot83}),
which allows us to identify
$\PSL$ with $ST\Hy^2$. A given orientation on $\Hy^2$ thus induces
an orientation on the $S^1$-fibres of $\PSL =ST\Hy^2\rightarrow\Hy^2$,
and hence on the $\R$-fibres of $\SL\rightarrow\Hy^2$. When we pass to
a left-quotient of $\SL$, these oriented $\R$-fibres descend to
oriented Seifert fibres. With this orientation convention, the
invariants of the multiple fibres are $(\alpha_j,1)$,
see~\cite[p.~467]{scot83}. On the other hand, there is a natural
{\em right\/} $S^1$-action on compact left-quotients of $\SL$. When this
right action is turned into a left action by the
inverse elements (while keeping the
orientation of $\SL$ and its quotient), the Seifert invariants
become $(\alpha_i,\alpha_i-1)$. This is the convention of
Raymond and Vasquez~\cite[pp.~169/70]{rava81}, which is
the more suitable one for our more algebraic considerations
in our earlier paper \cite{gego02} and below.
\end{rem}

A presentation of the fundamental group $\pi$ of $ST\Sigma$ is given by
\begin{eqnarray*}
\pi & = & \Bigl\{ u_1,v_1,\ldots ,u_g,v_g,q_1,\ldots ,q_n,h\co
\prod_i[u_i,v_i]\prod_jq_j=h^{2g-2},\\
&   & \mbox{}\;\;\;\;\;\;\;\; q_j^{\alpha_j}h^{\alpha_j-1}=1,\; h\;
\mbox{\rm central}\Bigr\}.
\end{eqnarray*}
Under the projection $ST\Sigma\rightarrow\Sigma$, the generators of
$\pi$ and $\pi^{\orb}$ correspond to each other as suggested by our
choice of notation. In other words, there is a representation
$\rho\in\calR (\pi ,\SL )$ with $\rho (\pi )=p^{-1}(\orho (\pi^{\orb}))$
and $p(\rho (u_i))=\orho (\ou_i))$ etc.
For further details see~\cite[Section~4]{gego02}.

The Seifert fibration $ST\Sigma\rightarrow\Sigma$, up to equivalence,
does not depend on the choice of hyperbolic structure on~$\Sigma$.
This allows us to speak of the unit tangent bundle $ST\Sigma$
(as a Seifert manifold) even when we have not fixed a metric on~$\Sigma$.
\section{Roots of the unit tangent bundle}
\label{section:roots}
Our aim is to classify $\rth$ roots of $ST\Sigma$ for $\Sigma$ an
oriented orbifold of hyperbolic type, by which we mean the
following.

\begin{defn}
An $\rth$ {\em root\/} of the unit tangent bundle $ST\Sigma$ is
an $r$-fold fibrewise connected and orientation preserving covering
$M\rightarrow ST\Sigma$ of $ST\Sigma$ by a Seifert manifold~$M$.
In other words, we require that each $S^1$-fibre
of $ST\Sigma$ is covered $r$ times positively by a single $S^1$-fibre of~$M$.
\end{defn}

\begin{rems}
(1) For $r=2$, such coverings are precisely the spin structures on~$\Sigma$.
Spin structures on orbifolds of arbitrary dimension were
defined and studied from the differential geometric point of view
(index theory, twistor theory) in \cite{dlm02}
and~\cite{bgr07}. The latter paper contains a general existence and
classification statement for spin structures on orbifolds, albeit only
for orbifolds whose singular set is of codimension at least~$4$.

(2) In the case of a principal $S^1$-bundle without multiple fibres,
one can pass to the associated complex line bundle. An $\rth$ root
then corresponds to a complex line bundle whose $\rth$ tensor
power is the given line bundle.
\end{rems}

Coverings $M\rightarrow ST\Sigma$ of the described type are automatically
regular, more specifically, $M$ is a principal $\Z_r$-bundle over~$ST\Sigma$.
The equivalence classes of (arbitrary) principal $\Z_r$-bundles
over $ST\Sigma$ are in natural one-to-one correspondence with the set
$\Hom (\pi ,\Z_r)$
of homomorphisms from the fundamental group $\pi$ of $ST\Sigma$ into~$\Z_r$,
the correspondence being given by associating with a principal $\Z_r$-bundle
its monodromy homomorphism~\cite[{\S}13.9]{stee51}. Thus, the fibrewise
connected and orientation preserving
coverings $M\rightarrow ST\Sigma$ are classified by the subset
\[ \Hom_1(\pi ,\Z_r):=\{ \delta\in\Hom (\pi ,\Z_r)\co\delta (h)=1\} .\]
If we drop the condition on orientations, we also have to
allow homomorphisms $\delta$ with $\delta (h)=-1$. This will become
relevant in Section~\ref{section:algebraic}.

\begin{rem}
In \cite[Remark~4.10]{gego02} we gave a detailed discussion of
the isomorphism between, on the one hand,
the deck transformation group of the universal covering $\widetilde{X}
\rightarrow X$ of a topological space $X$ and, on the other hand,
the fundamental group $\pi_1 (X,x_0)$. This isomorphism depends,
up to an inner automorphism, on the choice of a lift $\tilde{x}_0\in
\widetilde{X}$ of the base point~$x_0$. This dependence becomes irrelevant
once we consider homomorphisms into the abelian group~$\Z_r$. Thus, while
we usually think of $\pi$ as a deck transformation group, one may
still interpret the monodromy homomorphism $\pi\rightarrow\Z_r$
as being defined in terms of loops as in~\cite[{\S}13]{stee51}.
\end{rem}

We now want to give a characterisation of the homomorphisms
$\delta\in\Hom_1(\pi ,\Z_r)$ in terms of the allowable values on the
generators in the standard presentation of~$\pi$. In order to do so,
we need to recall a theorem of Raymond and Vasquez~\cite{rava81}
about the Seifert invariants of left-quotients of Lie groups,
cf.~\cite{gego95}. We have seen in the preceding section that,
once we equip $\Sigma$ with a hyperbolic structure, its unit tangent bundle
$ST\Sigma$ can be written as a left-quotient of $\SL$, and so the same
is true for its $r$-fold covering~$M$. 
Indeed, the fundamental group of the manifold $M$ corresponding to
$\delta\in\Hom_1(\pi ,\Z_r)$ is $\tpi =\ker\delta$. A representation
$\rho\in\calR (\pi ,\SL )$ as described at the end of
Section~\ref{section:tangent} induces a representation
$\trho\in\calR (\tpi ,\SL )$ of $\tpi$ as the deck transformation group
of~$M$.

By construction, $M$ is a Seifert manifold with $n$ multiple fibres
of multiplicities $\alpha_1,\ldots ,\alpha_n$ (just like the Seifert
manifold~$ST\Sigma$), but whereas the fibre index (see~\cite[Defn.~6]{gego95}
of $ST\Sigma$ equals~$1$, the fibre index of $M$ is~$r$. Then,
according to \cite{rava81} or~\cite{gego95}, the normalised Seifert invariants
\[ \{ g,b,(\alpha_1,\beta_1),\ldots ,(\alpha_n,\beta_n)\} \]
of $M$ (where $b$ is an integer and each $\beta_j$ an integer
between $1$ and $\alpha_j-1$)
are subject to the condition that there exist integers
$k_1,\ldots ,k_n$ such that
\begin{eqnarray}
rb       & = & 2g-2-\sum_{j=1}^n k_j,\label{eqn:rava1}\\
r\beta_j & = & \alpha_j-1+k_j\alpha_j,\;\; j=1,\ldots ,n.\label{eqn:rava2}
\end{eqnarray}
(Observe that these conditions are satisfied for $ST\Sigma$ with $r=1$,
$b=2g-2$, and all $k_j$ equal to zero.) For a given $\Sigma$, these
conditions impose severe restrictions on the possible values
of~$r$. These restrictions are implicit in~\cite{rava81}; for the
reader's convenience we deduce them directly from the equations
(\ref{eqn:rava1}) and~(\ref{eqn:rava2}).

\begin{lem}
\label{lem:rava}
If $r\in\N$ satisfies the Raymond--Vasquez relations (\ref{eqn:rava1})
and~(\ref{eqn:rava2}),
then $r$ is prime to $\alpha_1\cdots\alpha_n$ and divides
the integer $\alpha_1\cdots\alpha_n\cdot\chi^{\orb}$.

Conversely, if $r\in\N$ satisfies these latter conditions (for given
$g$, $n$ and $\alpha_j$), then there are integers $b$, $k_j$
and $\beta_j$ (with $1\leq\beta_j\leq\alpha_j-1$)
such that equations (\ref{eqn:rava1}) and (\ref{eqn:rava2}) are satisfied.
\end{lem}

\begin{proof}
From (\ref{eqn:rava2}) we see that $r$ must be prime to~$\alpha_j$. With
(\ref{eqn:rava1}) and (\ref{eqn:rava2}) one computes
\[ r\cdot\alpha_1\cdots\alpha_n\cdot
\bigl( b+\sum_{j=1}^n\frac{\beta_j}{\alpha_j}\bigr)
=-\alpha_1\cdots\alpha_n\cdot\chi^{\orb},\]
which proves the claimed divisibility.

For the converse, the condition $\gcd (r,\alpha_j)=1$ allows
us to choose integers $1\leq \beta_j\leq\alpha_j-1$ and $k_j$
such that (\ref{eqn:rava2}) holds. One then computes
\[ r\cdot\alpha_1\cdots\alpha_n\sum_{j=1}^n\frac{\beta_j}{\alpha_j}
=-\alpha_1\cdots\alpha_n\cdot\chi^{\orb}-
\alpha_1\cdots\alpha_n\cdot (2g-2-\sum_{j=1}^n k_j),\]
which shows that $r$ divides $2g-2-\sum_{j=1}^n k_j$, as was to be shown.
\end{proof}

\begin{rem}
Equation (\ref{eqn:rava2}) and the fact that $r$ and $\alpha_j$ are
coprime implies
that multiple fibres with the same $\alpha_j$ also have the same
$\beta_j$ (and hence the same~$k_j$). This is a unique feature
of left-quotients of $\SL$.
\end{rem}

The converse implication of the preceding lemma has the following
consequence. Given an $r\in\N$ satisfying the divisibility assumptions,
we find --- by the lemma --- a set of normalised Seifert
invariants satisfying the Raymond--Vasquez relations.
In particular, the Euler number
\[ e=-\bigl( b+\sum_{j=1}^n\frac{\beta_j}{\alpha_j}\bigr) \]
of the Seifert fibration must be non-zero, since $re=\chi^{\orb}<0$.
This means that the Seifert manifold $M$ defined by these invariants is
diffeomorphic to a
left-quotient of~$\SL$. The projection $\SL\rightarrow\PSL$
induces the Seifert fibration $M\rightarrow\Sigma$ over a hyperbolic
orbifold $\Sigma$ and gives
$M$ the structure of an $\rth$ root of $ST\Sigma$.

\begin{lem}
\label{lem:delta}
The homomorphisms $\delta\in\Hom_1(\pi,\Z_r)$ can take
arbitrary values on the generators $u_1,v_1,\ldots ,u_n,v_n$, but the value on
the $q_j$ is determined by $\delta (q_j)=k_j$ {\rm mod}~$r$.
\end{lem}

\begin{proof}
In $\Z_r$ we compute
\[ 0=\delta (1)=\delta (q_j^{\alpha_j}h^{\alpha_j-1})=
\delta (q_j)\alpha_j+\alpha_j-1. \]
From equation (\ref{eqn:rava2}) we see
that, first of all, $\alpha_j$ must be prime to~$r$, and secondly,
that $\delta (q_j)=k_j$ {\rm mod}~$r$, as claimed. Equation
(\ref{eqn:rava1}) implies that this condition on $\delta$ is consistent with
the other relation in the presentation of~$\pi$. It is then easy to
see that we may define a homomorphism $\delta\in\Hom_1(\pi,\Z_r)$
by prescribing arbitrary values on the $u_i$ and~$v_i$.
\end{proof}

We summarise our discussion in the following theorem.

\begin{thm}
\label{thm:existence}
The unit tangent bundle $ST\Sigma$ of an orbifold $\Sigma$ of
hyperbolic type admits an $\rth$ root if and only if $r\in\N$ is prime
to the multiplicities $\alpha_1,\ldots ,\alpha_n$ of the cone
points and a divisor of the integer
$\alpha_1\cdots\alpha_n\cdot\chi^{\orb}(\Sigma )$. In that case,
the distinct $\rth$ roots are in natural one-to-one correspondence with
the elements of $\Hom_1(\pi,\Z_r)$.
\hfill\qed
\end{thm}

\begin{rem}
There is a simple geometric explanation why $r$ needs to be
prime to the multiplicities $\alpha_1,\ldots ,\alpha_n$ for an
$\rth$ root of $ST\Sigma$ to exist. From the local model for
a fibre of multiplicity $\alpha_j$ one sees that when we pass
to an $r$-fold cover with connected covering of the multiple fibre,
then for $r$ not prime to $\alpha_j$ the covering of the regular fibres
will fail to be connected.
\end{rem}

Lemma~\ref{lem:delta} implies that any two homomorphisms in
$\Hom_1(\pi,\Z_r)$ differ by a homomorphism $\pi\rightarrow\Z_r$
that sends $h$ and the $q_j$ to zero. Such a homomorphism
may be interpreted as an element of
\[ \Hom (\pi_1(\Sigma ),\Z_r)\cong\Hom (H_1(\Sigma ),\Z_r)\cong
H^1(\Sigma;\Z_r)\cong\Z_r^{2g}.\]
(In the first term we really do mean the fundamental group $\pi_1(\Sigma)$,
not the orbifold fundamental group~$\pi^{\orb}$.)
This one-to-one correspondence of $\rth$ roots of $ST\Sigma$ with elements
of $H^1(\Sigma ;\Z_r)$, however, is not natural. All we have is a free
and transitive action of $H^1(\Sigma ;\Z_r)$ on the set of $\rth$ roots.

One way to give an explicit one-to-one correspondence between
$\Hom_1(\pi,\Z_r)$ and $\Z_r^{2g}$ is to fix a presentation for~$\pi$,
and then to associate with $\delta\in\Hom_1(\pi,\Z_r)$ the tuple
\[ (\delta (u_1),\delta (v_1),\ldots ,
\delta (u_g),\delta (v_g))\in\Z_r^{2g}.\]

\begin{rems}
(1) As observed by Johnson~\cite{john80}, there is a natural
geometric lifting of mod~$2$ homology classes from a surface to its unit
tangent bundle. Thus, spin structures on surfaces are {\em naturally\/}
classified both by $\Hom_1(\pi_1(\Sigma ),\Z_2)$ and $H_1(\Sigma;\Z_2)$.
There is no such natural lifting of mod~$r$ classes for $r$ greater
than~$2$. However, given a smooth simple closed curve on~$\Sigma$,
we can consider its tangential lift to $ST\Sigma$. This will be used
in the next section to help us understand the action of the
diffeomorphism group of $\Sigma$ on $\Hom_1(\pi ,\Z_r)$.

(2) For an honest $S^1$-bundle over an arbitrary manifold~$X$,
one can classify $\rth$ roots by mimicking the spectral
sequence argument of ~\cite[Chapter~II.1]{lami89} with $\Z_2$-coefficients
replaced by $\Z_r$-coefficients, cf.~\cite{kluk08}. This allows one to
show that an $\rth$ root exists if and only if the mod~$r$ reduction of the
Euler class of the $S^1$-bundle vanishes (which can also be seen by
more simple means), and then $\rth$ roots are
in (non-natural) one-to-one correspondence with the elements
of $H^1(X;\Z_r)$.
\end{rems}

We close this section by giving an explicit presentation of the
fundamental group $\tpi =\ker\delta$ of the manifold $M$ corresponding
to $\delta\in\Hom_1 (\pi ,\Z_r)$. Choose integers $s_i,t_i$
with $s_i\equiv \delta (u_i)$ and
$t_i\equiv\delta (v_i)$ mod~$r$, $i=1,\ldots ,g$. Then $\tpi$ is generated by
\[ \tu_i:=u_ih^{-s_i},\;
   \tv_i:=v_ih^{-t_i},\; i=1,\ldots ,g,\;\;
   \tq_j:=q_jh^{-k_j},\; j=1,\ldots ,n,\;\;
   \mbox{\rm and} \;\;
   \tih:=h^r.\]
With the help of the Raymond--Vasquez relations one sees that this yields
the presentation
\begin{eqnarray*}
\tpi & = & \Bigl\{ \tu_1,\tv_1,\ldots ,\tu_g,\tv_g,\tq_1,
\ldots ,\tq_n,\tih\co
\prod_i[\tu_i,\tv_i]\prod_j\tq_j=\tih^b,\\
  &   & \mbox{}\;\;\;\;\;\;\;\;\tq_j^{\alpha_j}\tih^{\beta_j}=1,\, \tih\;
\mbox{\rm central}\Bigr\}.
\end{eqnarray*}
\section{The action of diffeomorphisms on roots}
\label{section:diffeo}
We are now going to define an action of the diffeomorphism group of
$\Sigma$ on the set of $\rth$ roots of $ST\Sigma$. For $g\geq 2$, it will
be shown that this action is transitive for $r$ odd, and that it has exactly
two orbits for $r$ even; the case $g=1$ will
require an {\em ad hoc\/} treatment; on a given hyperbolic orbifold of
genus $g=0$ and for each $r$ there is at most one $\rth$ root,
since in that case $\Hom_1(\pi ,\Z_r)$ is trivial.
Throughout, we fix an orientation of
$\Sigma$, and the diffeomorphisms we consider are always
understood to be {\em orientation preserving}.
A diffeomorphism of an orbifold may at best permute cone points
of the same multiplicity. By Lemma~\ref{lem:delta} and
the remark following the proof of Lemma~\ref{lem:rava},
any such permutation can be
achieved by a diffeomorphism that induces the trivial action
on $\Hom_1(\pi ,\Z_r)$. So in order to understand
the action of the diffeomorphism
group of $\Sigma$ on $\Hom_1(\pi ,\Z_r)$, it suffices to
consider diffeomorphisms that fix a neighbourhood
of each cone point.

Let $M\rightarrow ST\Sigma$ be an $\rth$ root of $ST\Sigma$,
corresponding to some homomorphism $\delta\in\Hom_1(\pi ,\Z_r)$, and let $f$
be a diffeomorphism of $\Sigma$ as described. By slight abuse of notation,
we may regard the differential $Tf$ as a diffeomorphism of $ST\Sigma$;
the composition of the projection $M\rightarrow ST\Sigma$ with $Tf$
is then a new $\rth$ root of~$ST\Sigma$. We denote the corresponding
element in $\Hom_1(\pi ,\Z_r)$ by $f_*\delta$.

Geometrically this means the following. Given $u\in\pi$, represent it
by a loop in $ST\Sigma$. Then $(f_*\delta )(u)\in\Z_r$ is given by
the monodromy of the covering $M\rightarrow ST\Sigma$ along the
{\em preimage\/} of that loop under $Tf$.

We make one further abuse of notation. If $u$ is an oriented, smooth
closed curve on~$\Sigma$ (avoiding the cone points), we also write $u$ for its
tangential lift to a closed curve in $ST\Sigma$. Up to conjugation,
this represents a well-defined element in $\pi=\pi_1(ST\Sigma )$,
so it makes sense to speak of $\delta (u)\in\Z_r$. This abuse of notation
is justified by the fact that for $f$ a diffeomorphism of~$\Sigma$,
the tangential lift of $f(u)$ equals the image of the tangential lift
of $u$ under the differential~$Tf$.

\begin{figure}[h]
\labellist
\small\hair 2pt
\pinlabel $u_1$ [r] at 75 22
\pinlabel $u_2$ [r] at 218 22
\pinlabel $v_1$ [b] at 95 102
\pinlabel $v_2$ [b] at 239 102
\pinlabel $w_{12}$ [b] at 163 91
\endlabellist
\centering
\includegraphics[scale=0.7]{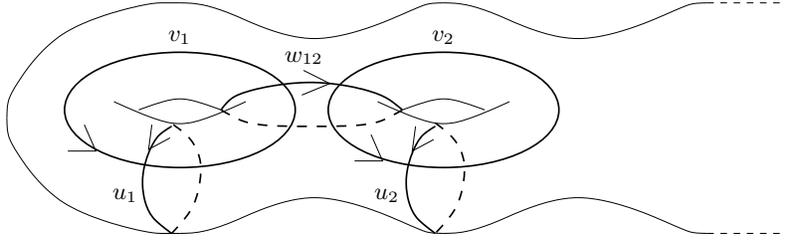}
  \caption{Loops on $\Sigma$ representing loops in $ST\Sigma$.}
  \label{figure:surface}
\end{figure}

Consider a topological model for $\Sigma$ as in
Figure~\ref{figure:surface}. Here $\Sigma$ is given the standard orientation,
so that the simple closed curves $u_i,v_i$ representing the
standard generators of $H_1(\Sigma)$ intersect positively in
a single point. The notation $u_i,v_i$ has been chosen in accordance with
the presentation of $\pi$ in Section~\ref{section:tangent}.
In the sequel we identify a homomorphism $\delta\in\Hom_1(\pi ,\Z_r)$
(representing an $\rth$ root of~$ST\Sigma$) with the corresponding
$2g$-tuple of integers mod~$r$, that is, we write
\[ \delta = (s_1,t_1,\ldots ,
s_g,t_g)\in\Z_r^{2g},\]
where $s_i=\delta (u_i)$ and $t_i=\delta (v_i)$.

For our discussion below we note that
the element $h\in\pi$ corresponds to a positively oriented regular fibre
of $ST\Sigma$, so it can be represented by a small positively oriented
circle on~$\Sigma$.

Next we want to show that for any given $\delta\in\Hom_1(\pi ,\Z_r)$
there is a diffeomorphism $f$ of $\Sigma$ such that $f_*\delta$
is in a very simple standard form. This is done by studying
the transformation behaviour of $\delta$ under certain Dehn twists
on~$\Sigma$.

For $u$ a {\em simple\/} closed curve on $\Sigma$, write
$f^u$ for the right-handed Dehn twist along~$u$.

\begin{lem}
\label{lem:transform}
Under the basic Dehn twists $f^{u_i}$, $f^{v_i}$ and $f^{w_{i,i+1}}$,
the tuple $\delta =(s_1,t_1,\ldots ,s_g,t_g)$ transforms as follows:
\begin{eqnarray*}
f^{u_i}_*\delta & = & (\ldots ,s_i,t_i-s_i,\ldots ),\\
f^{v_i}_*\delta & = & (\ldots , s_i+t_i,t_i,\ldots ),\\
f^{w_{i,i+1}}_*\delta & = & (\ldots , s_i,t_i-s_i+s_{i+1}-1,s_{i+1},
                             t_{i+1}+s_i-s_{i+1}+1,\ldots ).
\end{eqnarray*}
\end{lem}

\begin{proof}
The Dehn twist $f^{u_i}$ sends $u_i$ to itself and $v_i$ to
$u_i+v_i$; the differential $Tf^{u_i}$ has the same effect, when
those curves are regarded as loops in $ST\Sigma$. So the inverse
diffeomorphism sends $u_i$ to itself and $v_i$ to $v_i-u_i$.
This gives the formula for $f^{u_i}_*$. The argument for
$f^{v_i}_*$ is analogous.

In order to investigate $f^{w_{i,i+1}}_*$, we need to compute
$\delta (w_{i,i+1})$.

\begin{claim}
$\delta (w_{i,i+1})=\delta (u_i)-\delta (u_{i+1})+1.$
\end{claim}

This can be seen as follows (cf.~\cite{john80} for an analogous idea).
The disjoint union of the smooth curves $u_i$ and $-u_{i+1}$ (that is,
$u_{i+1}$ with reversed orientation) can be deformed smoothly into the
union of $w_{i,i+1}$ and a small circle oriented negatively (see
Figure~\ref{figure:u1u2} for a schematic illustration).
The tangential lift of the latter equals~$-h$. This implies
the claim.

\begin{figure}[h]
\labellist
\small\hair 2pt
\pinlabel $u_i$ [l] at 2 20
\pinlabel $-u_{i+1}$ [r] at 74 126
\pinlabel $w_{i,i+1}$ [b] at 614 92
\pinlabel $-h$ [l] at 628 73
\endlabellist
\centering
\includegraphics[scale=0.5]{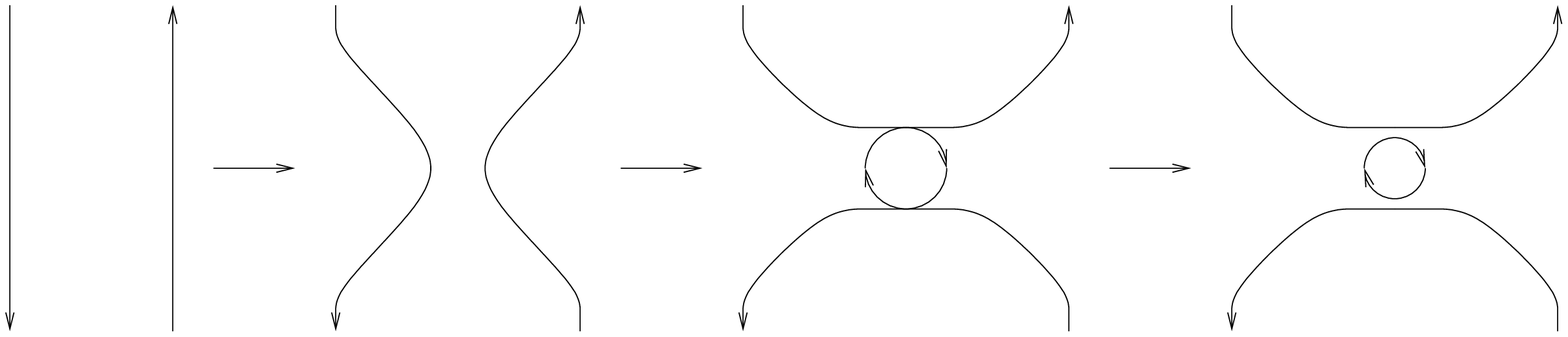}
  \caption{The smooth deformation from $u_i-u_{i+1}$ to $w_{i,i+1}-h$.}
  \label{figure:u1u2}
\end{figure}

Now, the inverse of $f^{w_{i,i+1}}$ sends $v_i$ to $v_i-w_{i,i+1}$,
and $v_{i+1}$ to $v_{i+1}+w_{i,i+1}$; the other basic loops remain
unchanged. In conjunction with the claim, this gives the formula
for~$f^{w_{i,i+1}}_*$.
\end{proof}

\begin{rem}
The formul\ae\ of Lemma~\ref{lem:transform} --- for the case $r=2$, where
signs do not matter --- where derived earlier by D\c{a}browski
and Percacci~\cite{dape86} by quite involved calculations
in local coordinates. Related considerations can also be found in the
work of Sipe~\cite{sipe86}. She studied $\rth$ roots of the unit tangent
bundle of hyperbolic {\em surfaces\/} with the aim of describing certain
finite quotients of their mapping class group.
\end{rem}

The signs in the formul\ae\ of Lemma~\ref{lem:transform}
change when we perform left-handed Dehn twists.
Therefore, Dehn twists along $u_i$ and $v_i$
enable us to perform Euclid's algorithm on any pair of integers
representing the pair $(s_i,t_i)$ of mod~$r$ classes.
This implies that we can reduce one component to zero and
the other to the unique element $d_i\in\Z_r$ determined by the
conditions that the principal ideal in $\Z_r$ generated by $d_i$
equal the ideal generated by $s_i$ and~$t_i$,
and that the integer representative of $d_i$ lying between $1$
and $r$ be a divisor of~$r$. By slight abuse of notation we write
this last condition as $d_i|r$.
The pair $(d_i,0)$ can be changed to the pair $(0,d_i)$ by
further such Dehn twists.
In total, we can find a composition of
Dehn twists of $\Sigma$ that transforms $\delta$ to
\[ (0,d_1,\ldots ,0,d_g).\]

In order to simplify this further, we have to bring the curves
$w_{i,i+1}$ into play. By the claim, the transformed $\delta$
takes the value $1$ on $w_{i,i+1}$. Thus, when we perform
$d_1$ right-handed Dehn twist along $w_{12}$, the tuple
$(0,d_1,0,d_2,\ldots )$ changes to $(0,0,0,d_1+d_2,\ldots )$.
Continuing with the appropriate Dehn twists along $w_{23}$ up
to $w_{g-1,g}$, we find a diffeomorphism transforming $\delta$ to
\[ (0,\ldots ,0,d_1+\cdots +d_g).\]
We shall presently describe further
Dehn twists that bring $\delta$ into one of the forms listed in
the next proposition.

\begin{prop}
\label{prop:standard}
By a sequence of Dehn twists, we can bring $\delta$ into one
of the following standard forms:
\begin{itemize}
\item $(0,\ldots ,0,0)$ if $g\geq 2$ and $r$ odd,
\item $(0,\ldots ,0,0)$ or $(0,\ldots ,0,1)$ if $g\geq 2$ and $r$ even,
\item $(0,d)$ with $d|r$ if $g=1$ (beware that this includes $d=
r\equiv 0$).
\end{itemize}
\end{prop}

Of course, for $g=1$ the surface $\Sigma$
will be of hyperbolic type only if there is at least one cone point. For
$g=0$ (and at least three cone points), $\Hom_1(\pi ,\Z_r)$ is trivial.

\begin{proof}[Proof of Proposition~\ref{prop:standard}]
The case $g=1$ has been settled by the discussion preceding the
proposition. In the case $g\geq 2$,
 we may assume that $\delta$ has already been
transformed into the form $(0,\ldots ,0,d)$, as yet without any information
on~$d$.

We only write the last four
components of the $2g$-tuple in $\Z_r^{2g}$. We claim that there
are Dehn twists giving the following sequence of transformations:
\[ (0,0,0,d)\longrightarrow (0,\mp 1 ,0, d\pm 1)\longrightarrow
(0,\pm 1,0,d\pm 1)\longrightarrow (0,0,0,d\pm 2).\]
Indeed, the first and third step are given by a Dehn twist (of the
appropriate sign) along $w_{g-1,g}$, the second by a sequence of Dehn twists
along $u_{g-1}$ and $v_{g-1}$.

So we can always reduce the last component to $0$ or~$1$. If $r$ is odd,
then Dehn twists along $u_g$ and $v_g$ allow us to transform from
$(0,1)$ (in the last two components) to $(0,2)$ --- since either of $1$
or $2$ generates the same principal ideal in~$\Z_r$, namely the
full ring.
\end{proof}

The standard forms listed in the preceding proposition turn out to
be pairwise inequivalent under the action of the
diffeomorphism group. We first show this for the case $g=1$.

\begin{lem}
\label{lem:standard-1}
For $g=1$, two standard forms $(0,d)$ and $(0,d')$, where we
think of $d,d'$ as integers between $1$ and $r$ (which divide~$r$), are
diffeomorphic if and only if $d=d'$.
\end{lem}

\begin{proof}
Assume without loss of generality that $d'\leq d$.
The action of the diffeomorphism group of a
hyperbolic orbifold $\Sigma$
of genus $1$ translates into the standard $\mbox{\rm SL}(2,\Z )$-action
on $\Z_r^2=\Hom_1(\pi ,\Z_r)$. The orbit of $(0,d)$ under
this action consists of elements of the form $(md,nd)$ (with $m$ and $n$
coprime). Since $d$ is a divisor of~$r$, the number $nd$ (thought
of as an integer) can be congruent to $d'$ mod~$r$ only if
$d$ is a divisor of~$d'$, which forces $d=d'$.
\end{proof}

The $\Z_2$-invariant that distinguishes the standard forms in the case
$g\geq 2$ (and $r$ even) goes back to Atiyah~\cite{atiy71}. A spin structure
on an honest surface $\Sigma$ has an associated complex line bundle~$L$.
Once a complex structure has been chosen on~$\Sigma$, one can speak
of holomorphic sections of~$L$. The dimension mod~$2$
of the vector space of holomorphic sections turns out to be
independent of the chosen complex structure; this is Atiyah's
invariant of spin structures. As remarked earlier, Johnson~\cite{john80}
defined a natural lifting of mod~$2$ homology classes from a surface $\Sigma$
to its unit tangent bundle. A spin structure on $\Sigma$ then
gives rise to a quadratic form on $H_1(\Sigma;\Z_2)$. Johnson goes
on to show that the Arf invariant of that quadratic form (whose
definition can be found on any Turkish 10 Lira note) equals
Atiyah's invariant.

\begin{rem}
The $2$-dimensional spin cobordism group
$\Omega_2^{\mbox{\scriptsize\rm spin}}$ is isomorphic to $\Z_2$;
the Atiyah invariant distinguishes the two cobordism classes.
\end{rem}

Now we allow once again arbitrary orbifolds $\Sigma$ of hyperbolic type.
Motivated by Johnson's work, we define a $\Z_2$-valued invariant of
an $\rth$ root $\delta$ of $ST\Sigma$ (with $r$ even), which we write
as $\delta = (s_1,t_1,\ldots ,s_g,t_g)\in\Z_r^{2g}$, by
\[ A(\delta )=\sum_{i=1}^g (s_i+1)(t_i+1)\;\;\mbox{\rm mod}\; 2.\]
Note that, for $r$ even, this mod~$2$ reduction is well defined.
The definition of this invariant can also be phrased as follows.
Given a principal $\Z_r$-bundle $M\rightarrow ST\Sigma$ with $r$ even,
there is an intermediate double covering of~$ST\Sigma$. Thus, an $\rth$ root
(with $r$ even) induces in a natural way a spin structure. The
$A$-invariant is simply the Atiyah invariant of that spin structure.

\begin{lem}
The number $A(\delta )\in\Z_2$ is a diffeomorphism invariant,
i.e.\ for any (orientation preserving) diffeomorphism $f$ of $\Sigma$ one
has $A(f_*\delta )=A(\delta )$.
\end{lem}

\begin{proof}
We need only consider diffeomorphisms that fix the cone points.
The group of such diffeomorphisms is generated by the Dehn twists
along $u_i,v_i$ and $w_{i,i+1}$. The invariance of $A(\delta )$
under these Dehn twists can be checked easily with the formul\ae\
in Lemma~\ref{lem:transform}.
\end{proof}

Obviously, the two standard forms $(0,\ldots ,0,0)$ and
$(0,\ldots ,0,1)$ (for $g\geq 2$ and $r$ even) are distinguished
by the $A$-invariant.

\begin{defn}
For $g\geq 2$ and $r$ even, we say an $\rth$ root $\delta$ is
of {\em even\/} (resp.\ {\em odd\/}) {\em type\/} if $A(\delta )$ equals~$0$
(resp.~$1$).
\end{defn}

So the standard form $(0,\ldots ,0,0)$ is of even type for $g$ even,
and of odd type for $g$ odd; the standard form $(0,\ldots ,0,1)$
has the complementary type.

\begin{prop}
\label{prop:number}
For $g\geq 2$ and $r$ even, the number of $\rth$ roots of
even (resp.\ odd) type equals $r^{2g}(2^g\pm 1)/2^{g+1}$.
\end{prop}

\begin{proof}
Write $r=2s$. An $\rth$ root $\delta =(s_1,t_1,\ldots ,s_g,t_g)$
will be even if and only if an even number of summands $(s_i+1)(t_i+1)$
in $A(\delta )$ are odd. Such a summand is odd if and only if
both $s_i$ and $t_i$ are even, which gives us $s^2$ possibilities for
choosing $s_i$ and~$t_i$. On the other hand, there are
$3s^2$ possibilities for choosing $s_i$ and $t_i$ such that
$(s_i+1)(t_i+1)$ becomes even. It follows that the number
of roots of even type is given by
\begin{eqnarray*}
\sum_{k\;\mathrm{even}}\binom{g}{k}
(s^2)^k(3s^2)^{g-k}
 & = & \frac{1}{2}\bigl( (s^2+3s^2)^g+(-s^2+3s^2)^g\bigr)\\
 & = & \frac{1}{2}\bigl( (4s^2)^g+(2s^2)^g\bigr)\\
 & = & r^{2g}(2^g+1)/2^{g+1}.
\end{eqnarray*}
For roots of odd type, the calculation is analogous.
\end{proof}

\begin{rem}
In the case $r=2$, i.e.\ for spin structures,
Propositions \ref{prop:standard} and \ref{prop:number} are
well known --- especially, it seems, among mathematical physicists.
Our arguments for deriving
them generalise those of D\c{a}browski and Percacci~\cite{dape86}.
An alternative approach can be found in the work of
Alvarez-Gaum\'e, Moore and Vafa~\cite{amv86}. They appeal to the
relation between spin structures and theta functions in order
to describe the action of the diffeomorphism group.
\end{rem}
\section{An algebraic reformulation}
\label{section:algebraic}
The Baer--Nielsen theorem for the orbifold $\Sigma$ says, in essence,
that the group of all (not just orientation preserving)
diffeomorphisms of $\Sigma$ modulo those isotopic to
the identity can be identified with the group $\Out (\pi^{\orb})$
of outer automorphisms of~$\pi^{\orb}$; see~\cite{zies73}.
We now want to use this to reformulate the action of the diffeomorphism
group on the space of $\rth$ roots of $ST\Sigma$ in an algebraic way.
This serves as a preparation for the next section, where we tie
up our discussion of $\rth$ roots with the moduli problem
for so-called taut contact circles, which was addressed
in our earlier paper~\cite{gego02}. As announced there, the
results of the present note allow us to count the connected
components of the moduli spaces in question.

There is an obvious action of $\Aut (\pi )$ on
\[ \Hom_{\pm 1}(\pi ,\Z_r):= \{ \delta\in\Hom (\pi ,\Z_r)
\co \delta (h)=\pm 1\}.\]
(The fact that we now allow $\delta (h)=-1$ corresponds to
having orientation reversing diffeomorphisms included
in the discussion.)
This descends to an action of $\Out (\pi )$, since $\Z_r$ is
abelian. Thus, in order to define an action of $\Out (\pi^{\orb})$
on $\Hom_{\pm 1}(\pi ,\Z_r)$, we should first define a suitable
lift from $\Aut (\pi^{\orb})$ to $\Aut (\pi )$.
Recall from \cite[Lemma~4.13]{gego02} that there is a short exact
sequence
\[ 0\longrightarrow\Z^{2g}\longrightarrow\Aut (\pi )
\longrightarrow\Aut (\pi^{\orb})\longrightarrow 1.\]
(This holds true for the fundamental group $\pi$ of any Seifert
manifold which is a left quotient of $\SL$ and has base orbifold~$\Sigma$.)
Thus, algebraically, it is not clear how to define a lifting.
Instead, we find a suitable lift by a direct appeal to the
Baer--Nielsen theorem. Put briefly, we represent a given
element of $\Out (\pi^{\orb})$ by an orbifold diffeomorphism $f$
of~$\Sigma$, and then find the lift as the automorphism corresponding to
the differential~$Tf$. From that construction it is clear that
our algebraic definition of the action by the diffeomorphism group on the
set of $\rth$ roots corresponds to the geometric definition in the
preceding sections (except that we have replaced a left action by a right
action, which is owed to the conventions in the algebraic setting
of the next section).

\begin{lem}
\label{lem:algebraic}
There is a natural right action of\/ $\Out (\pi^{\orb})$
on\/ $\Hom_{\pm 1}(\pi ,\Z_r)$, defined as follows. Given
a class $[\otheta ]\in\Out (\pi^{\orb})$, represented by an
automorphism $\otheta\in\Aut (\pi^{\orb})$, there is a geometrically
defined lifting of this representative to an automorphism
$\vartheta\in\Aut (\pi )$. Then the action of $[\otheta ]$
on $\delta\in\Hom_{\pm 1}(\pi ,\Z_r)$ is defined by $\delta\mapsto
\delta\circ\vartheta$.
\end{lem}

\begin{proof}
By the Nielsen theorem~\cite[Theorem~8.1]{zies73} there is a
diffeomorphism $f$ of $\Sigma$ (fixing a base point~$x_0$)
covered by a diffeomorphism
$\tf$ of $\tSigma$ (fixing a chosen lift $\widetilde{x}_0$
of~$x_0$) such that
\[ \tf\circ \ou\circ \tf^{-1}=\otheta (\ou )\;\;\;\mbox{\rm for all}\;\;
\ou\in\pi^{\orb}.\]
Regard the differential $Tf$ as a diffeomorphism of $ST\Sigma$,
and let $\widetilde{Tf}$ be a lift to a diffeomorphism
of~$\SL$. Define $\vartheta\in\Aut (\pi )$ by
\[ \vartheta (u)= \widetilde{Tf}\circ u\circ\widetilde{Tf}^{-1}
\;\;\;\mbox{\rm for all}\;\; u\in\pi .\]

Since the fibre class $h$ generates the centre of~$\pi$, we have
$\vartheta (h)=h^{\pm 1}$. So the homomorphism $\delta\circ\vartheta$
is still an element of $\Hom_{\pm 1}(\pi ,\Z_r)$. Moreover, the definitions
imply that $\otheta_1\circ\otheta_2$ lifts to $\vartheta_1\circ
\vartheta_2$,
so the prescription $\delta\mapsto\delta\circ\vartheta$ does indeed
define a {\em right\/} action, provided we can establish independence of
choices.

Two different lifts of $Tf$ differ by a deck transformation
of $ST\Sigma$, i.e.\ an element of $\pi$. So the corresponding
lifts $\vartheta$ differ by an inner automorphism of~$\pi$.
Thus, the homomorphism $\delta\circ\vartheta$ into the abelian group $\Z_r$ is
independent of this choice of lift.

Next, we show that $\delta\circ\vartheta$ depends only on the
class~$[\otheta ]$, not on the choice of representative~$\otheta$,
or in other words, that for any {\em inner\/}
automorphism $\otheta$ we have $\delta\circ\vartheta
=\delta$. By the Baer theorem~\cite[Theorem 3.1]{zies73}, the
Nielsen realisation $f$ of any inner automorphism $\otheta$ of $\pi^{\orb}$
is isotopic to the identity (by an isotopy {\em not\/} fixing
the base point, in general). Then $Tf$ is likewise isotopic to
the identity. This isotopy lifts to a fibre isotopy between
$\widetilde{Tf}$ and a deck transformation of $\SL\rightarrow ST\Sigma$.
This implies that the resulting $\vartheta$ will be an inner
automorphism of~$\pi$, and hence $\delta\circ\vartheta =\delta$, as we wanted
to show.

Finally, it remains to verify that the construction does not depend on the
choice of Nielsen realisation $f$. Two such realisations differ by
a diffeomorphism whose lift to $\tSigma$ induces the identity
on~$\pi^{\orb}$. Then the argument concludes as before by an appeal to
Baer's theorem.
\end{proof}

\section{The moduli space of taut contact circles}
\label{section:moduli}
Let $M$ be a given closed, orientable $3$-manifold diffeomorphic
to a left quotient of $\SL$ with fundamental group~$\tpi$.
This is in a unique way a Seifert manifold over an orbifold
$\Sigma$ of hyperbolic type, with a well-defined fibre index~$r$.
Recall from the end of Section~\ref{section:roots}
the presentation of $\tpi$ involving the normalised Seifert
invariants of~$M$.

As shown in our paper~\cite{gego02}, the Teichm\"uller space
$\T (M)$ of taut contact circles, i.e.\ the space of taut contact circles
on $M$ modulo diffeomorphisms isotopic to the identity,
can be identified with $\Inn (\SL)\backslash\calR (\tpi ,\SL )$,
where $\calR$ stands for the Weil space of representations
as in Section~\ref{section:orbifold}.
The moduli space $\M (M)$ of taut contact circles, i.e.\ the space of taut
contact circles on $M$ modulo all diffeomorphisms of~$M$,
is in turn given by $\T (M)/\Out (\tpi )$. With this algebraic translation
taken for granted, nothing further needs to be known about
taut contact circles (not even their definition), i.e.\ the following
can be read as a discussion of these algebraically defined spaces,
where we want to understand the action of $\Out (\tpi )$
on $\T (M)=\Inn (\SL )\backslash\calR (\tpi ,\SL )$ with the
help of the geometry of $\rth$ roots of $ST\Sigma$.
See also ~\cite{klr85} for the relevance of such questions to the deformation
theory of Seifert manifolds.

\begin{rem}
To a large extent we follow the notational conventions of~\cite{gego02}.
The one difference that needs to be pointed out is that in our
previous paper, $\pi$ denoted the fundamental group of~$M$,
as $ST\Sigma$ did not play much of a role in our discussion there.
In the present paper, $\pi$ denotes the fundamental group of $ST\Sigma$,
and $\tpi$ that of~$M$.
\end{rem}

Write $\T (\Sigma )$ for the Teichm\"uller space of hyperbolic metrics
on the base orbifold~$\Sigma$, together with a choice of
orientation. This means that $\T (\Sigma )$ has two connected
components $\T^+(\Sigma )$ and $\T^-(\Sigma )$. Algebraically,
$\T (\Sigma )$ may be thought of as $\Inn (\PSL )\backslash\calR
(\pi^{\orb},\PSL )$.
In Section~4 of \cite{gego02} it was shown that
$\T (M)$ is a trivial principal $\Z^{2g}$-bundle over~$\T (\Sigma )$.
For $\Aut (\tpi )$ there is a short exact sequence as for
$\Aut (\pi )$ in the previous section.
The normal subgroup $\Z^{2g}\subset\Aut (\tpi )$ acts as
$(r\Z )^{2g}$ on the mentioned principal bundle. This implies
that $\T (M)/\Z^{2g}$ --- where the quotient is taken
under the action of $\Z^{2g}\subset\Aut (\tpi )$ ---
is a trivial $r^{2g}$-fold covering
of $\T (\Sigma )$, and the moduli space of taut contact circles
on $M$ can be described as
\[ \M (M)=\bigl(\T (M)/\Z^{2g}\bigr) /\Out (\pi^{\orb}).\]
So the following theorem essentially settles the moduli problem for taut
contact circles on left quotients of~$\SL$. Here $\pi$ denotes, as before,
the fundamental group of $ST\Sigma$. For the proof below, notice
that there are quotient maps $\tpi\rightarrow\pi^{\orb}$ and
$\pi\rightarrow\pi^{\orb}$, given by quotienting out the normal subgroup
generated by the central element $\tih$ and~$h$, respectively.

\begin{thm}
\label{thm:teich}
The quotient $\T (M)/\Z^{2g}$ of the Teichm\"uller space
of taut contact circles on $M$ under the action of
$\Z^{2g}\subset\Aut (\tpi )$ has a natural description
as follows:
\[ \T (M)/\Z^{2g} = \Hom_1(\pi ,\Z_r)\times\T^+(\Sigma ) \sqcup
\Hom_{-1}(\pi ,\Z_r)\times\T^-(\Sigma )\]
On the second factors $\T^{\pm}(M)$, the right action of $\Out (\pi^{\orb})$
is the obvious one; on the first factors $\Hom_{\pm 1}(\pi,\Z_r)$,
the group $\Out (\pi^{\orb})$ acts from
the right as described in Section~\ref{section:algebraic}.
\end{thm}

\begin{rem}
If the Nielsen realisation
$f$ of an automorphism $\otheta$ of $\pi^{\orb}$ is orientation reversing
(so that $\otheta$ will exchange the components $\T^{\pm}(\Sigma )$),
then the differential~$Tf$, regarded as a diffeomorphism of $ST\Sigma$,
will reverse the fibre direction,
so $\otheta$ will also exchange $\Hom_{\pm 1}(\pi ,\Z_r)$.
In fact, no left-quotient of $\SL$ admits any orientation
reversing diffeomorphism~\cite{nera78}.
\end{rem}

\begin{proof}[Proof of Theorem~\ref{thm:teich}]
First we are going to define a map from the left-hand side
$\T (M)/\Z^{2g}$ to the right factors $\T^+(\Sigma )
\sqcup\T^-(\Sigma )=\T (\Sigma )$ on the right-hand side.
Recall from \cite[Section~4]{gego02} that the projection
$\SL\rightarrow\PSL$ induces a covering map $\calR (\tpi ,\SL )
\rightarrow\calR (\pi^{\orb},\PSL )$, which in turn induces
a well-defined map $\T (M)\rightarrow \T (\Sigma )$, since
any inner automorphism of $\SL$ induces an inner automorphism
of $\PSL$. The action of $\Z^{2g}\subset\Aut (\tpi )$ on $\tpi$
is given by multiplying the generators $\tu_i,\tv_i$ with
the corresponding power of the central element~$\tih$.
Since this central element generates the kernel of the
quotient map $\tpi\rightarrow\pi^{\orb}$, we get an induced map
$\T (M)/\Z^{2g}\rightarrow\T (\Sigma )$.

Next we want to define a map $\T (M)/\Z^{2g}\rightarrow\Hom_{\pm 1}
(\pi,\Z_r)$ to the left factors on the right-hand side.
This means that, given $\trho\in\calR (\tpi ,\SL )$ representing an element
$[\trho ]\in\T (M)/\Z^{2g}$, and given~$u\in\pi$,
we want to define $\delta (u)\in
\Z_r$ in such a way that $\delta$ becomes a homomorphism $\pi\rightarrow\Z_r$
sending $h$ to~$\pm 1$, and such that $\delta$ is independent of
the chosen representative~$\trho$.

Thus, start with $\trho$ and $u$ as described. The element $u\in\pi$
projects to an element $\ou\in\pi^{\orb}$, which in turn lifts
to an element $\tu\in\tpi$, unique up to powers of~$\tih$. Likewise,
the representation $\trho\in\calR (\tpi ,\SL )$ projects to
a representation
\[ \orho\in\calR (\pi^{\orb},\PSL )=\calR^+(\pi^{\orb},\PSL )
\sqcup\calR^-(\pi^{\orb},\PSL ), \]
as observed in the
first part of the proof, and then can be lifted in a preferred way
to a representation $\rho\in\calR (\pi,\SL )$;
cf.\ Section~\ref{section:orbifold} for the notation~$\calR^{\pm}$.

In \cite[Section~4]{gego02} we gave a definition of such a
preferred lift that also allowed us to lift from a
representation of $\pi^{\orb}$ to one of~$\tpi$. Here, where we only
want to lift to a representation of~$\pi$, we shall make a choice
that leads to a natural description of the $\Out (\pi^{\orb})$-action.
In order to allow unique lifting of maps to universal covers,
we choose base points in $\Sigma$, $ST\Sigma$ and their universal covers
in such a way that all relevant projections are
base point preserving. Likewise, we choose a base point
in $\SL=\widetilde{ST\Hy}$ over a base point in~$\Hy$; this
determines a base point in any discrete quotient of these spaces.

Now to the definition of~$\rho$. In the sequel it is understood that
all diffeomorphisms are base point preserving.
Choose a diffeomorphism $g\co\Sigma\rightarrow \orho (\pi^{\orb})
\backslash\Hy^2$ whose (unique) lift $\tg$ to the universal
cover satisfies
\[ \tg\circ \ou\circ \tg^{-1}=\orho (\ou )\;\;\;\mbox{\rm for all}\;\;
\ou\in\pi^{\orb}.\]
This is possible by the Nielsen theorem again; observe the formal similarity
with the argument in the proof of Lemma~\ref{lem:algebraic}.
Now, with $L$ denoting left multiplication in $\SL$,
define the preferred lift $\rho$ of $\orho$ by
\[ L_{\rho (u)}= \widetilde{Tg}\circ u\circ\widetilde{Tg}^{-1}
\;\;\;\mbox{\rm for all}\;\; u\in\pi .\]

\begin{rem}
The preferred lift as defined in \cite{gego02} depended on
a choice of presentation of~$\pi$. If we take the $u_i$ and $v_i$
as the tangential lifts of $\ou_i$ and $\ov_i$, then the
preferred lift defined here is the same as that in~\cite{gego02}.
\end{rem}

When we identify $\Hy$ with the upper half-plane in $\C$,
and $\SL$ with $\Hy^2\times\R$ with coordinates $(z,\theta )$,
cf.~\cite[p.~58]{gego02},
we can describe the left action of $\rho (u)$ on $\SL$ explicitly
(at least for some elements $u\in\pi$).
For $h$ there is no choice in the lifting; one has
\[ \rho (h)(z,\theta )=(z,\theta\pm 2\pi ),\]
where the sign is determined by $\orho\in\calR^{\pm}(\pi^{\orb},\PSL )$.
Similarly, one has
\[ \trho (\tih )(z,\theta )=(z,\theta\pm 2\pi r).\]
The lift $\rho (q_j)$ is completely determined by the relation which
$q_j$ satisfies in the group~$\pi$. For $u_i$ resp.\ $v_i$,
any lift other than the preferred one $\rho (u_i)$ resp.\
$\rho (v_i)$ would differ from it
by an arbitrary translation
in the $\theta$-component by integer multiples of~$2\pi$.
Moreover, the action of $\bw\in\Z^{2g}\subset\Aut (\tpi )$ on
$\trho$ is given by $\trho\mapsto\trho_{r\bw}$, with
\begin{eqnarray*}
\trho_{r\bw}(\tu_i)(z,\theta ) & = & \trho (\tu_i)(z,\theta )+
                                    (0,2\pi rw_{2i-1}),\\
\trho_{r\bw}(\tv_i)(z,\theta ) & = & \trho (\tv_i)(z,\theta )+
                                    (0,2\pi rw_{2i}).
\end{eqnarray*}

Now back to the construction of the homomorphism~$\delta$ corresponding
to the class $[\trho ]\in\T (M)/\Z^{2g}$. Since both $\rho (u)$ and
$\trho (\tu )$ are lifts of $\orho (\ou )\in\PSL$ to
$\SL$, their action on the $\theta$-component differs
by a shift by some integer multiple of~$2\pi$, so we can define
$\delta (u)\in\Z$ by
\begin{equation}
\label{eqn:delta}
\rho (u)(z,\theta )= \trho (\tu )(z,\theta )+ (0,2\pi\delta (u)).
\end{equation}
Since $\rho$ is fixed to be the preferred lift of~$\orho$, the only
ambiguity in this equation is the lift $\tu$ of~$\ou$, which may
be changed by powers of~$\tih$. From the described action of $\trho (\tih )$
we conclude that $\delta (u)$ is well defined mod~$r$, so we may regard it
as a map into~$\Z_r$. By construction it is clear that $\delta$ has the
homomorphism property. For $u=h$ we may choose
$\tu =1$; this gives $\delta (h)=\pm 1$, where the sign again
corresponds to $\orho\in\calR^{\pm}$, as it should. Hence
$\delta\in\Hom_{\pm 1}(\pi,\Z_r)$.

Inner automorphisms of $\SL$ act trivially on the $\theta$-component,
so $\delta$ only depends on the class of $\trho$ in $\T (M)$.
Moreover, $\delta (u)$ does not change mod~$r$ when
$\trho$ is replaced by some $\trho_{r\bw}$ in the same
orbit under the $\Z^{2g}$-action on $\T (M)$. This finishes the
construction of the map
\[ \T (M)/\Z^{2g} \rightarrow \Hom_1(\pi ,\Z_r)\times\T^+(\Sigma ) \sqcup
\Hom_{-1}(\pi ,\Z_r)\times\T^-(\Sigma ).\]

We show this map to be a bijection by exhibiting an explicit
inverse. The defining equation~(\ref{eqn:delta}) for $\delta$
can be read backwards, as it were, in order to define the desired
inverse map. Thus, given $\orho\in\calR^{\pm}(\pi^{\orb},\PSL )$
and $\delta\in\Hom_{\pm 1}(\pi,\Z_r)$ (with matching signs),
we would like to use (\ref{eqn:delta}) to define $\trho$.
This is indeed possible, if we take a little care. First of
all, we know that there is no choice in defining
$\trho (\tih )$ and $\trho (\tq_j)$, so we only need to
consider elements $\tu\in\tpi$ which are not stabilised under
the $\Z^{2g}$-action $\trho (\tu )\mapsto \trho_{r\bw}(\tu )$.
Let $\ou\in\pi^{\orb}$ be the projection of~$\tu$, and $u\in\pi$
a lift of~$\ou$. In the equation
\[ \trho (\tu )(z,\theta )=\rho (u) (z,\theta )-(0,2\pi\delta (u)), \]
with $\rho$ taken as the preferred lift of~$\orho$,
the right-hand side can be made sense of if the $\theta$-component
is read as lying in $\R/2\pi r\Z$, and it does not depend on the
choice of lift~$u$. So for $\tu$ of the described kind, we can use this
equation (given $\orho$, $\delta$ and~$\tu$)
to get a well-defined element $[\trho ]\in\T (M)/\Z^{2g}$.
This prescription obviously defines an inverse of the previously
constructed map.

It remains to show that the right action of $\Out (\pi^{\orb})$ is
as claimed in the theorem. Given $[\otheta ]\in\Out (\pi^{\orb})$,
let $\ttheta\in\Aut (\tpi )$ be any lift of $\otheta$,
and $\vartheta\in\Aut (\pi )$ the lift constructed in
the proof of Lemma~\ref{lem:algebraic}.
The action of $[\otheta ]$ on $\T (M)/\Z^{2g}$
is given by $\trho\mapsto\trho\circ\ttheta$. This is indeed
well defined: the choice of representative $\otheta$ of the
class $[\otheta ]$ is irrelevant, because in $\T (M)$
we have taken the quotient under $\Inn (\SL )$;
the specific lifting to $\ttheta$ is of no importance
in the quotient $\T (M)/\Z^{2g}$. That the action
of $\Out (\pi^{\orb})$ on the right-hand side of the identity
in the theorem is also as claimed now follows from equation~(\ref{eqn:delta})
and the observation that our construction of the preferred lift
of $\orho$ entails that $\rho\circ\vartheta$ is the
preferred lift of~$\orho\circ\otheta$.

This concludes the proof of Theorem~\ref{thm:teich}.
\end{proof}

\begin{rem}
With $\rho$ being the preferred lift of~$\orho$,
equation~(\ref{eqn:delta}) is precisely the algebraic reformulation
of the geometric definition of $\delta$ as a monodromy
homomorphism given in Section~\ref{section:roots}.
\end{rem}

When we take the quotient under the action of
$\Out (\pi^{\orb})$, the trivial covering $\T (M)/\Z^{2g}\rightarrow
\T^+(\Sigma )\sqcup\T^-(\Sigma )$ given by Theorem~\ref{thm:teich}
becomes a possibly branched
covering $\M (M)\rightarrow\M (\Sigma )$, where
$\M (\Sigma )=\T (\Sigma )/\Out (\pi^{\orb})$ denotes the moduli
space of hyperbolic metrics on~$\Sigma$. 

We are now interested in the number of connected
components of $\M (M)$, and the number of sheets in each
connected component of the covering $\M (M)\rightarrow \M (\Sigma )$.
The space
$\M (\Sigma )$ is connected, so the number of
connected components of $\M (M)$ equals the number of orbits
of the $\Out (\pi^{\orb})$-action on $\Hom_{\pm 1} (\pi ,\Z_r)$.
Geometrically, this corresponds to the number of orbits
of the action on $\Hom_1(\pi ,\Z_r)$ given by the
orientation preserving diffeomorphisms of~$\Sigma$.
Moreover, the number of sheets in each connected component
of the covering $\M (M)\rightarrow \M (\Sigma )$
is given by the length of the corresponding orbit.

So the following theorem, the larger part of which
was announced in \cite{gego02},
is a direct consequence of
Propositions \ref{prop:standard} and~\ref{prop:number},
and Lemma~\ref{lem:standard-1}.
(As before, we write $r$ for the fibre index of the unique
Seifert fibration $M\rightarrow\Sigma$; the genus of $\Sigma$
is denoted by~$g$.)

\begin{thm}
\label{thm:moduli}
The moduli space $\M (M)$ of taut contact circles on $M$
is a branched covering over the moduli space $\M (\Sigma )$
of hyperbolic metrics on~$\Sigma$.

For $g=0$, the covering map $\M (M)\rightarrow\M (\Sigma )$
is a homeomorphism.

For $g=1$, the number of connected components of $\M (M)$
equals the number of divisors of $r$.
The number of sheets in the component of $\M (M)$ corresponding
to $d|r$ equals the number of ordered pairs $(s,t)$ of integers
mod~$r$ that generate the same ideal in $\Z_r$ as~$d$.

For $g\geq 2$ and $r$ odd, $\M (M)$ is connected, and the branched
covering $\M (M)\rightarrow\M (\Sigma )$ has $r^{2g}$ sheets.

For $g\geq 2$ and $r$ even, $\M (M)$ has two connected components,
and the number of sheets in the two components equals
$r^{2g}(2^g\pm 1)/2^{g+1}$.
\qed
\end{thm}


\begin{thebibliography}{99}
%
\bibitem{amv86}
{L. Alvarez-Gaum\'e, G. Moore and C. Vafa}:
{\it Theta functions, modular invariance, and strings},
{Comm. Math. Phys.}
{\bf 106} (1986), 1--40.
%
\bibitem{atiy71}
{M. F. Atiyah}:
{\it Riemann surfaces and spin structures},
{Ann. Sci. \'Ecole Norm. Sup. (4)}
{\bf 4} (1971), 47--62.
%
\bibitem{bgr07}
{F. A. Belgun, N. Ginoux and H.-B. Rademacher}:
{\it A singularity theorem for twistor spinors},
{Ann. Inst. Fourier (Grenoble)}
{\bf 57} (2007), 1135--1159.
%
\bibitem{dape86}
{L. D\c{a}browski and R. Percacci}:
{\it Spinors and diffeomorphisms},
{Comm. Math. Phys.}
{\bf 106} (1986), 691--704.
%
\bibitem{dlm02}
{C. Dong, K. Liu and X. Ma}:
{\it On orbifold elliptic genus};
in {Orbifolds in Mathematics and Physics} (Madison, 2001),
Contemp. Math. {\bf 310},
American Mathematical Society, Providence, 2002, 87--105.
%
\bibitem{gego95a}
{H. Geiges and J. Gonzalo P\'erez}:
{\it Contact geometry and complex surfaces},
{Invent. Math.}
{\bf 121} (1995), 147--209.
%
\bibitem{gego95}
{H. Geiges and J. Gonzalo P\'erez}:
{\it Seifert invariants of left-quotients of $3$-dimensional simple
Lie groups},
{Topology Appl.}
{\bf 66} (1995), 117--127.
%
\bibitem{gego02}
{H. Geiges and J. Gonzalo P\'erez}:
{\it Moduli of contact circles},
{J. Reine Angew. Math.}
{\bf 551} (2002), 41--85.
%
\bibitem{jarv00}
{T. J. Jarvis}:
{\it Geometry of the moduli of higher spin curves},
{Internat. J. Math.}
{\bf 11} (2000), 637--663.
%
\bibitem{john80}
{D. Johnson}:
{\it Spin structures and quadratic forms on surfaces},
{J. London Math. Soc. (2)}
{\bf 22} (1980), 365--373.
%
\bibitem{kluk08}
{M. Klukas}:
{Engelstrukturen},
Diplomarbeit,
Universit\"at zu K\"oln, 2008.
%
\bibitem{klr85}
{R. Kulkarni, K. B. Lee and F. Raymond}:
{\it Deformation spaces for Seifert manifolds};
in {Geometry and Topology} (College Park, 1983/84),
Lecture Notes in Math. {\bf 1167}, Springer-Verlag,
Berlin, 1985, 180--216.
%
\bibitem{lami89}
{H. B. Lawson, Jr. and M.-L. Michelsohn}:
{Spin Geometry},
Princeton Math. Ser. {\bf 38},
Princeton University Press, Princeton, 1989.
%
\bibitem{nera78}
{W. D. Neumann and F. Raymond}:
{\it Seifert manifolds, plumbing, $\mu$-invariant and orientation
reversing maps};
in {Algebraic and Geometric Topology} (Santa Barbara, 1977),
Lecture Notes in Math. {\bf 664},
Springer-Verlag, Berlin, 1978, 163--196.
%
\bibitem{poli04}
{A. Polishchuk}:
{\it Witten's top Chern class on the moduli space of higher
spin curves};
in {Frobenius Manifolds} (Bonn, 2002),
Aspects Math. {\bf E36},
Vieweg, Wiesbaden, 2004, 253--264.
%
\bibitem{rava81}
{F. Raymond and A. T. Vasquez}:
{\it $3$-manifolds whose universal coverings are Lie groups},
{Topology Appl.}
{\bf 12} (1981), 161--179.
%
\bibitem{scot83}
{P. Scott}:
{\it The geometries of $3$-manifolds},
{Bull. London Math. Soc.}
{\bf 15} (1983), 401--487.
%
\bibitem{sipe86}
{P. L. Sipe}:
{\it Some finite quotients of the mapping class group of a surface},
{Proc. Amer. Math. Soc.}
{\bf 97} (1986), 515--524.
%
\bibitem{stee51}
{N. Steenrod}:
{The Topology of Fibre Bundles},
Princeton Math. Ser. {\bf 14},
Princeton University Press, Princeton, 1951.
%
\bibitem{witt93}
{E. Witten}:
{\it Algebraic geometry associated with matrix models of two dimensional
gravity};
in {Topological Methods in Modern Mathematics}
(Stony Brook, 1991),
Publish or Perish, Houston, 1993, 235--269.
%
\bibitem{zies73}
{H. Zieschang}:
{\it On the homeotopy group of surfaces},
{Math. Ann.}
{\bf 206} (1973), 1--21.
%
\end{thebibliography}
\end{document}